\documentclass{amsart}

\def\C{\mathbb C}
\def\R{\mathbb R}
\def\X{\mathbb X}

\def\Z{\mathbb Z}

\def\Y{\mathbb Y}
\def\Z{\mathbb Z}
\def\N{\mathbb N}

\newtheorem{theorem}{Theorem}[section]
\newtheorem{lemma}[theorem]{Lemma}

\theoremstyle{definition}
\newtheorem{definition}[theorem]{Definition}
\newtheorem{example}[theorem]{Example}

\theoremstyle{remark}
\newtheorem{remark}[theorem]{Remark}

\numberwithin{equation}{section}

\begin{document}

\title[Pseudo Almost Periodic Solutions]{Existence of Pseudo Almost Periodic Solutions \\to Some
Classes of Partial Hyperbolic\\ Evolution Equations}

\author{Toka Diagana}
\address{Department of Mathematics, Howard University, 2441 6th Street N.W., Washington, DC 20059, USA}

\email{tdiagana@howard.edu}





\begin{abstract}
The paper examines the existence of pseudo almost periodic
solutions to some classes of partial hyperbolic evolution
equations. Namely, sufficient conditions for the existence and
uniqueness of pseudo almost periodic solutions to those classes of
hyperbolic evolution equations are given. Applications include the
existence of pseudo almost periodic solutions to the transport and
heat equations with delay.
\end{abstract}

\maketitle

\section{Introduction}
Let $(\X , \|\cdot\|)$ be a Banach space and let $A: D(A) \subset
\X \mapsto \X$ be a sectorial linear operator (Definition
\ref{sect}). For $\alpha \in (0 , 1)$, the space $\X_\alpha$
denotes an abstract {\it intermediate} Banach space between $D(A)$
and $\X$. Examples of those $\X_\alpha$ include, among others, the
fractional spaces $D((-A)^\alpha)$ for $\alpha \in (0 ,1)$, the
reel interpolation spaces $D_A(\alpha,\infty)$ due to J. L. Lions
and J. Peetre, and the H\"older spaces $D_A(\alpha),$ which
coincide with the continuous interpolation spaces that had been
introduced in the literature by G. Da Prato and P. Grisvard.

In \cite{cm,dmn, DMNT,lhl},  some sufficient conditions for the
existence and uniqueness of pseudo almost periodic solutions to
the abstract (semilinear) differential equations,
\begin{eqnarray}\label{1}
u'(t) + A u(t) = f(t , u(t)), \ \ t \in \R, \ \ \mbox{and}
\end{eqnarray}
\begin{eqnarray}\label{11111}
u'(t) + A u(t) = f(t , Bu(t)), \ \ t \in \R,
\end{eqnarray}
where $-A$ is a Hille-Yosida linear operator (respectively, the
infinitesimal generator of an analytic semigroup, and the
infinitesimal generator of a $C_0$-semigroup), $B$ is a densely
defined closed linear operator on $\X$, and $f: \R \times \X
\mapsto \X$ is a jointly continuous function, were given.
Similarly, in \cite{DT}, some sufficient conditions for the
existence and uniqueness of pseudo almost periodic solutions to
the class of partial evolution equations
\begin{eqnarray}\label{0}
\displaystyle{\frac{d}{dt}} \left[u(t) + f(t, Bu(t))\right] = A
u(t) + g(t, Cu(t)),\; \; t \in \R
\end{eqnarray}
where $A$ is the infinitesimal generator of an exponentially
stable semigroup acting on $\X$, $B, C$ are arbitrary densely
defined closed linear operators on $\X$, and $f, g$ are some
jointly continuous functions, were given.

Now note the assumptions made in \cite{DT} require much more
regularity for the operator $A$, that is, being the infinitesimal
generator of an analytic semigroup. In this paper we address such
an issue by studying pseudo almost periodic solutions to (\ref{0})
in the case when $A$ is a sectorial operator whose corresponding
analytic semigroup $(T(t))_{t \geq 0}$ is hyperbolic,
equivalently,
$$\sigma(A) \cap i \R = \emptyset.$$ Note that (\ref{0}) in the case when $A$ is sectorial corresponds to
several interesting situations encountered in the literature.
Applications include, among others, the existence and uniqueness
of pseudo almost periodic solutions to the hyperbolic transport
and heat equations with delay.

As in \cite{BMN, GT} in this paper we consider a general
intermediate space $\X_\alpha$ between $D(A)$ and $\X$. In
contrast with the fractional power spaces considered in some
recent papers of the author et al. \cite{dmn, DMNT}, the
interpolation and H\"older spaces, for instance, depend only on
$D(A)$ and $\X$ and can be explicitly expressed in many concrete
cases. The literature related to those intermediate spaces is very
extensive, in particular, we refer the reader to the excellent
book by A. Lunardi \cite{lunardi}, which contains a comprehensive
presentation on this topic and related issues.

The concept of pseudo almost periodicity, which is the central
question in this paper was introduced in the literature in the
early nineties by C. Zhang \cite{z1, z2, z3} as a natural
generalization of the well-known Bohr almost periodicity. Thus
this new concept is welcome to implement another existing
generalization of almost periodicity, that is, the concept of
asymptotically almost periodicity due to Fr\'echet~\cite{cor,
fink}.

The existence of almost periodic, asymptotically almost periodic,
and pseudo almost periodic solutions is one of the most attractive
topics in qualitative theory of differential equations due to
their significance and applications in physics, mathematical
biology, control theory, physics and others.

Some contributions on almost periodic, asymptotically almost
periodic, and pseudo almost periodic solutions to abstract
differential and partial differential equations have recently been
made in \cite{aea1, aea2, am, cm, dia, dmn, DMNT, DT, lhl}.
However, the existence and uniqueness of pseudo almost periodic
solutions to (\ref{0}) in the case when $A$ is sectorial is an
important topic with some interesting applications, which is still
an untreated question and this is the main motivation of the
present paper. Among other things, we will make extensive use of
the method of analytic semigroups associated with sectorial
operators and the Banach's fixed-point principle to derive
sufficient conditions for the existence and uniqueness of a pseudo
almost periodic (mild) solution to \eqref{0}.

\section{Preliminaries}
This section is devoted to some preliminary facts needed in the
sequel. Throughout the rest of this paper, $(\X, \|\cdot\|)$
stands for a Banach space, $A$ is a sectorial linear operator
(Definition~\ref{sect}), which is not necessarily densely defined,
and $B, C$ are (unbounded) linear operators such that $A + B + C$
is not trivial as each solution to \eqref{0} belongs to
$D(A+B+C)$. Now if $A$ is a linear operator on $\X$, then
$\rho(A)$, $\sigma(A)$, $D(A)$, $N(A)$, $R(A)$ stand for the
resolvent, spectrum, domain, kernel, and range of $A$. The space
$B(\X, \Y)$ denotes the Banach space of all bounded linear
operators from $\X$ into $\Y$ equipped with its natural norm with
$B(\X, \X) = B(\X)$.

\subsection{Sectorial Linear Operators and their Associated Semigroups}
\begin{definition}\label{sect}
A linear operator $A: D(A) \subset \X \mapsto \X$ (not necessarily
densely defined) is said to be sectorial if the following hold:
there exist constants $\omega\in \R$, $\displaystyle \theta\in
(\frac{\pi}{2},\pi)$, and $M>0$ such that
\begin{align}\label{sect}
&\qquad\rho(A)\supset S_{\theta,\omega}:=\{\lambda\in \C:\lambda\neq\omega,\ \ |\arg(\lambda-\omega)|<\theta\}, \ \ \mbox{and}\\
&\qquad \|R(\lambda,A)\|\leq \frac{M}{|\lambda-\omega|},  \quad
\lambda\in S_{\theta,\omega}.
\end{align}
\end{definition}

The class of sectorial operators is very rich and contains most of
classical operators encountered in the literature. Two examples of
sectorial operators are given as follows:

\begin{example}
Let $p \geq 1$ and let $\X = L^p(\R)$ be the Lebesgue space
equipped with its norm $\|\cdot\|_p$ defined by
$$\|\varphi\|_p = \left(\int_{\R} |\varphi(x)|^p dx\right)^{1/p}.$$

Define the linear operator $A$ on $L^p(\R)$ by
$$D(A) = W^{2,p}(\R), \ \ A(\varphi) = \varphi^{\prime\prime}, \ \ \forall \varphi \in D(A).$$

It can be checked that the operator $A$ is sectorial on $L^p(\R)$.
\end{example}

\begin{example}
Let $p \geq 1$ and let $\Omega \subset \R^d$ be open bounded
subset with $C^2$ boundary $\partial \Omega$. Let $\X:=
L^p(\Omega)$ be the Lebesgue space equipped with the norm,
$\|\cdot\|_p$ defined by,
$$\|\varphi\|_p = \left(\int_{\Omega} |\varphi(x)|^p dx\right)^{1/p}.$$

Define the operator $A$ as follows:
$$D(A)= W^{2, p}(\Omega) \cap W^{1,p}_{0}(\Omega), \ \ A(\varphi)= \Delta \varphi, \ \ \forall \varphi \in D(A),$$
where $\displaystyle \Delta = \sum_{k=1}^d
\frac{\partial^2}{\partial x_k^2}$ is the Laplace operator.

It can be checked that the operator $A$ is sectorial on
$L^p(\Omega)$.

\end{example}

It is well-known that \cite{lunardi} if $A$ is sectorial, then it
generates an analytic semigroup $(T(t))_{t\geq0}$, which maps
$(0,\infty)$ into $B(\X)$ and such that there exist $M_0, M_1 > 0$
with
\begin{align}
&\|T(t)\|\leq M_0 e^{\omega t}, \quad t> 0,\\
&\|t(A-\omega)T(t)\|\leq M_1 e^{\omega t}, \quad t>
0.\label{analy}
\end{align}

Throughout the rest of the paper, we suppose that the semigroup
$(T(t))_{t\geq0}$ is hyperbolic, that is, there exist a projection
$P$ and constants $M, \delta>0$ such that $T(t)$ commutes with
$P$, $N(P)$ is invariant with respect to $T(t)$, $T(t):
R(Q)\mapsto R(Q)$ is invertible, and the following hold
\begin{equation}\label{hyP}
\|T(t)Px\|\leq Me^{-\delta t}\|x\| \qquad \mbox{ for } t\geq0,
\end{equation}
 \begin{equation}\label{hyQ}
  \|T(t)Qx\|\leq
Me^{\delta t}\|x\| \qquad \mbox{ for } t\leq0, \end{equation}
 where $Q:=I-P$ and,  for $ t\leq0$,  $T(t):= (T(-t))^{-1} $.

Recall that the analytic semigroup $(T(t))_{t\geq0}$ associated
with $A$ is hyperbolic if and only if $$\sigma(A)\cap
i\R=\emptyset,$$ see, e.g., \cite[Prop. 1.15, pp.305]{nag}.

\begin{definition}
Let $\alpha\in(0,1)$. A Banach space $(\X_\alpha,
\|\cdot\|_\alpha)$ is said to be an intermediate
 space between $D(A)$ and $\X$, or a space of class ${\mathcal J}_\alpha$,
 if $D(A)\subset \X_\alpha\subset \X$ and there is a constant  $c>0$
 such that
 \begin{equation}\label{extrap}
 \|x\|_{\alpha}\leq
c\|x\|^{1-\alpha}\|x\|_{A}^{\alpha},\qquad x\in D(A),
\end{equation}
where $\|\cdot\|_A$ is the graph norm of $A$.
\end{definition}
Concrete examples of $\X_\alpha$ include $D((-A^\alpha))$ for
$\alpha\in(0,1)$, the domains
 of the the fractional powers of $A$,  the real interpolation
 spaces $D_A(\alpha,\infty)$,  $\alpha\in(0,1)$,  defined as follows
 $$\left\{%
\begin{array}{ll}
\displaystyle    D_A(\alpha,\infty):=\{x\in \X: [x]_\alpha=\sup_{0< t\leq1}\|t^{1-\alpha}AT(t)x\|<\infty\}\\
\displaystyle    \|x\|_{\alpha}=\|x\|+[x]_{\alpha},
\end{array}%
\right. $$ the abstract H\"older spaces
$D_A(\alpha):=\overline{D(A)}^{\|.\|_\alpha}$ as well as the
complex interpolation spaces $[\X,D(A)]_\alpha$, see A. Lunardi
\cite{lunardi} for details.

  For a hyperbolic analytic  semigroup $(T(t))_{t\geq0}$, one can easily check that
  similar estimations as both \eqref{hyP} and \eqref{hyQ} still hold with norms
  $\|\cdot\|_\alpha$. In fact,  as the part of $A$ in $R(Q)$ is
  bounded, it follows from \eqref{hyQ} that
  \begin{equation*}
  \|AT(t)Qx\|\leq
C'e^{\delta t}\|x\| \qquad \mbox{ for } t\leq0. \end{equation*}

 Hence, from \eqref{extrap} there exists a constant $c(\alpha)>0$  such that
   \begin{equation}\label{hyQ*}
  \|T(t)Qx\|_{\alpha}\leq
c(\alpha)e^{\delta t}\|x\| \qquad \mbox{ for } t\leq0.
\end{equation}

In addition to the above, the following holds
  \begin{equation*}
\|T(t)Px\|_{\alpha}\leq \|T(1)\|_{B(\X,\X_\alpha)}\|T(t-1)Px\|
\qquad \mbox{ for } t\geq1,
\end{equation*}
 and hence from \eqref{hyP}, one obtains
\begin{equation*}
\|T(t)Px\|_{\alpha}\leq M'e^{-\delta t}\|x\|,  \qquad
t\geq1,\end{equation*} where $M'$ depends on $\alpha$. For
$t\in(0,1]$, by \eqref{analy} and \eqref{extrap}
$$ \|T(t)Px\|_{\alpha}\leq M'' t^{-\alpha}\|x\|. $$

Hence, there exist constants $M(\alpha)>0$ and $\gamma>0$ such
that
  \begin{equation}\label{hyP*}
\|T(t)Px\|_{\alpha}\leq M(\alpha)t^{-\alpha}e^{-\gamma t}\|x\|
\qquad \mbox{ for } t>0.
\end{equation}

\subsection{Pseudo Almost Periodic Functions}
Let $(\Y, \|\cdot\|_\Y)$ be another Banach space. Let $BC(\R ,
\X)$ (respectively, $BC(\R \times \Y, \X)$) denote the collection
of all $\X$-valued bounded continuous functions (respectively, the
class of jointly bounded continuous functions $F: \R \times \Y
\mapsto \X$).  The space $BC(\R, \X)$ equipped with its natural
norm, that is, the sup norm defined by
$$\displaystyle \|u\|_\infty = \sup_{t \in \R} \|u(t)\|$$ is a
Banach space. Furthermore, $C(\R, \Y)$ (respectively, $C(\R \times
\Y, \X)$) denotes the class of continuous functions from $\R$ into
$\Y$ (respectively, the class of jointly continuous functions $F:
\R \times \Y \mapsto \X$).

\begin{definition}\label{D}
A function $f \in C(\R , \X)$ is called (Bohr) almost periodic if
for each $\varepsilon > 0$ there exists $l(\varepsilon)
> 0$ such that every interval of length  $l(\varepsilon)
$ contains a number $\tau$ with the property that
$$\|f(t +\tau)
- f(t) \| < \varepsilon \ \ \mbox{for each} \ \ t \in \R.$$
\end{definition}

The number $\tau$ above is called an {\it
$\varepsilon$-translation} number of $f$, and the collection of
all such functions will be denoted $AP(\X)$.

\begin{definition}\label{D}
A function $F \in C(\R \times \Y, \X)$ is called (Bohr) almost
periodic in $t \in \R$ uniformly in $y \in \Y$ if for each
$\varepsilon
> 0$ and any compact $K \subset \Y$ there exists $l(\varepsilon)$ such that every interval of length  $l(\varepsilon)$ contains
a number $\tau$ with the property that
$$\|F(t + \tau, y) - F(t, y)\| < \varepsilon \ \ \mbox{for each} \ \ t \in \R, \ \ y \in K.$$

The collection of those functions is denoted by $AP(\R \times
\Y)$.
\end{definition}

Set

$$AP_0(\X) := \{ f \in BC(\R , \X): \lim_{r \to \infty}
\displaystyle{\frac{1}{2r}} \int_{-r}^r \| f(s)\| ds = 0\},$$ and
define $AP_0(\R \times \X)$ as the collection of functions $ F \in
BC(\R \times \Y, \X)$ such that
$$\lim_{r \to \infty}
\displaystyle{\frac{1}{2r}} \int_{-r}^r \| F(t,u)\| dt = 0$$
uniformly in $u \in \Y.$

\begin{definition}\label{DDD}
A function $f \in BC(\R , \X)$ is called pseudo almost periodic if
it can be expressed as $f = g + \phi,$ where $g \in AP(\X)$ and
$\phi \in AP_0(\X)$. The collection of such functions will be
denoted by $PAP({\mathbb X})$.
\end{definition}

\begin{remark}
The functions $g$ and $\phi$ in Definition \ref{DDD} are
respectively called the {\it almost periodic} and the {\it ergodic
perturbation} components of $f$. Moreover, the decomposition given
in Definition \ref{DDD} is unique.
\end{remark}

Similarly,

\begin{definition}
A function $F \in C(\R \times \Y, \X)$ is said to pseudo almost
periodic in $t \in \R$ uniformly in $y \in \Y$ if it can be
expressed as $F = G + \Phi,$ where $G \in AP(\R \times \Y)$ and
$\phi \in AP_0(\R \times \Y)$. The collection of such functions
will be denoted by $PAP(\R \times {\mathbb \Y})$.
\end{definition}

\section{Main results}
To study the existence and uniqueness of pseudo almost periodic
solutions to (\ref{0}) we need to introduce the notion of mild
solution to it.

\begin{definition} Let $\alpha \in (0 , 1)$.
A bounded continuous function $u: \R \mapsto \X_\alpha$ is said to
be a mild solution to (\ref{0}) provided that the function $s\to
AT(t-s)Pf(s,Bu(s))$ is integrable on $(-\infty,t)$, $s\to
AT(t-s)Qf(s,Bu(s))$ is integrable on $(t, \infty)$
 for each $t\in \mathbb{R},$ and
\begin{eqnarray}
u(t)&=&-f(t,Bu(t)) - \int_{-\infty}^{t}AT(t-s)P f(s, Bu(s))ds
\nonumber\\
&+& \int_{t}^{\infty}AT(t-s)Qf(s, Bu(s))ds +
\int_{-\infty}^{t}T(t-s)Pg(s, Cu(s))ds \nonumber\\
&-& \int_{t}^{\infty}T(t-s)Qg(s, Cu(s))ds\nonumber
\end{eqnarray}
for each $\forall t \in \R$.

\end{definition}

Throughout the rest of the paper we denote by $\Gamma_1, \Gamma_2,
\Gamma_3,$ and $\Gamma_4$, the nonlinear integral operators
defined by
$$(\Gamma_1 u)(t) := \int_{-\infty}^{t}AT(t-s)P f(s, Bu(s))ds,
\ \ (\Gamma_2 u)(t) := \int_{t}^{\infty}AT(t-s)Qf(s, Bu(s))ds,$$
$$(\Gamma_3 u)(t) :=\int_{-\infty}^{t}T(t-s)Pg(s, Cu(s))ds, \ \
\mbox{and}$$ $$ (\Gamma_4 u)(t) :=\int_{t}^{\infty}T(t-s)Qg(s,
Cu(s))ds .$$

To study \eqref{0} we require the following assumptions:

\begin{itemize}
\item[({\bf H1})] The operator $A$ is sectorial and generates a
hyperbolic (analytic) semigroup $(T(t))_{t\geq0}$.

\item[({\bf H2})] Let $0<\alpha<1$. Then $\X_\alpha=
D((-A^\alpha))$, or $\X_\alpha= D_A(\alpha,p), 1\leq p\leq
+\infty$, or $\X_\alpha= D_A(\alpha)$, or
$\X_\alpha=[\X,D(A)]_\alpha$. Let $B,C:\X_\alpha \longrightarrow
\X$ be bounded linear operators.

\item[({\bf H3})] Let $1>\beta>\alpha$, and $f: \R\times
\X\longrightarrow \X_{\beta}$ be a pseudo almost periodic function
in $t \in \R$ uniformly in $u \in \X$, $g: \R \times \X \mapsto
\X$ be pseudo almost periodic in $t \in \R$ uniformly in $u \in
\X$.

\item[({\bf H4})] The functions $f, g$ are uniformly
    Lipschitz with respect to the second argument in the following
    sense: there exists $K > 0$ such that
    $$
    \|f(t,u)-f(t,v)\|_\beta\leq K \|u-v\|,
    $$
 and
$$\|g(t,u)-g(t,v)\|\leq K \|u-v\|
    $$
    for all $u,v\in \X$ and $t\in \R$.

\end{itemize}
In order to show that $\Gamma_1$ and $\Gamma_2$ are well defined,
we need the following estimates.

\begin{lemma}\label{beta} Let $0<\alpha,\beta<1$. Then
\begin{align}
  \|AT(t)Qx\|_\alpha&\leq
ce^{\delta t}\|x\|_\beta \qquad \mbox{ for }
t\leq0,\label{beta1}\\
 \|AT(t)Px\|_{\alpha}&\leq
ct^{\beta-\alpha-1}e^{-\gamma t}\|x\|_\beta, \qquad \mbox{ for }
t>0.\label{beta2}
\end{align}
\end{lemma}

\begin{proof}
As for \eqref{hyQ*}, the fact that the part of $A$ in $R(Q)$ is
  bounded yields
  \begin{equation*}
  \|AT(t)Qx\|\leq
ce^{\delta t}\|x\|_\beta,  \quad   \|A^2T(t)Qx\|\leq ce^{\delta
t}\|x\|_\beta,\qquad \mbox{ for } t\leq0,
\end{equation*}
since $X_\beta\hookrightarrow X.$
 Hence, from \eqref{extrap} there is a constant $c(\alpha)>0$  such that
   \begin{equation*}
  \|AT(t)Qx\|_{\alpha}\leq
c(\alpha)e^{\delta t}\|x\|_\beta \qquad \mbox{ for } t\leq0.
\end{equation*}
Furthermore,
  \begin{align*}
\|AT(t)Px\|_{\alpha}&\leq
\|AT(1)\|_{B(\X,\X_\alpha)}\|T(t-1)Px\|\\
&\leq ce^{-\delta t}\|x\|_\beta,   \qquad \mbox{ for } t\geq1.
\end{align*}

Now for $t\in(0,1]$, by \eqref{analy} and \eqref{extrap}, one has
$$ \|AT(t)Px\|_{\alpha}\leq c t^{-\alpha-1}\|x\|,$$
and
$$ \|AT(t)Px\|_{\alpha}\leq c t^{-\alpha}\|Ax\|,$$
for each $x\in D(A)$. Thus, by reiteration Theorem (see
\cite{lunardi}), it follows that
$$ \|AT(t)Px\|_{\alpha}\leq c t^{\beta-\alpha-1}\|x\|_\beta$$
for every $x\in \X_\beta$ and $0<\beta<1$, and hence, there exist
constants $M(\alpha)>0$ and $\gamma>0$ such that
  \begin{equation*}
\|T(t)Px\|_{\alpha}\leq M(\alpha)t^{\beta-\alpha-1}e^{-\gamma
t}\|x\|_\beta \qquad \mbox{ for } t>0.
\end{equation*}
\end{proof}

\begin{lemma}\label{ll2} Under assumptions {\rm ({\bf H1})-({\bf H2})-({\bf H3})-({\bf
H4})}, the integral operators $\Gamma_3$ and $\Gamma_4$ defined
above map $PAP(\X_\alpha)$ into itself.
\end{lemma}

\proof Let $u \in PAP(\X_\alpha)$. Since  $C \in B(\X_\alpha, \X)$
it follows that $Cu \in PAP(\X)$. Setting  $h(t) = g(t, Cu(t))$
and using the theorem of composition of pseudo almost periodic
functions \cite[Theorem 5]{am} it follows that $h \in PAP(\X)$.
Now, write $h = \phi + \zeta$ where $\phi \in AP(\X)$ and $\zeta
\in AP_0(\X)$. Thus $\Gamma_3u$ can be rewritten as
\begin{eqnarray*}
(\Gamma_3u)(t) &=&\int_{-\infty }^tT(t-s)P\phi(s)ds +
\int_{-\infty }^tT(t-s)P\zeta(s)ds.
\end{eqnarray*}

Set $$\displaystyle \Phi(t) = \int_{-\infty }^tT(t-s)P\phi(s)ds,$$
and $$\displaystyle \Psi(t)= \int_{-\infty }^tT(t-s)P\zeta(s)ds$$
for each $t \in \R$. The next step consists of showing that $\Phi
\in AP(\X_\alpha)$ and $\Psi \in AP_0(\X_\alpha)$.

Clearly, $\Phi \in AP(\X_\alpha)$. Indeed, since $\phi \in
AP(\X)$, for every $\varepsilon > 0$ there exists $l(\varepsilon)
> 0$ such that for all $\xi$ there is $\tau \in [ \xi , \xi
+ l(\varepsilon)]$ with
$$\|\Phi(t +\tau) - \Phi(t) \| < \mu
\ . \ \varepsilon \ \ \mbox{for each} \ t \in \R,$$ where $\mu =
\displaystyle{\frac{\gamma^{1-\alpha}}{M(\alpha) \Gamma
(1-\alpha)}}$ with $\Gamma$ being the classical gamma function.

Now using the expression $$\displaystyle \Phi(t + \tau) - \Phi(t)
= \int_{-\infty}^{t} T(t -s) P\left(\phi(s + \tau) -
\phi(s)\right) ds$$ and \eqref{hyP*} it easily follows that $$\|
\Phi(t + \tau) - \Phi(t) \|_\alpha < \varepsilon \ \ \mbox{for
each} \ t \in \R,$$and hence, $\Phi \in AP(\X_\alpha)$. To
complete the proof for $\Gamma_3$, we have to show that $t \mapsto
\Psi(t)$ is in $AP_0(\X_\alpha)$. First, note that $s \mapsto
\Psi(s)$ is a bounded continuous function. It remains to show that

\begin{align*}
\lim_{r\to \infty} \displaystyle{\frac{1}{2r}} \ \int_{-r}^r \|
\Psi(t) \|_\alpha \, dt = 0.
\end{align*}

 Again using \eqref{hyP*} one obtains that
\begin{align*}
\lim_{r \to \infty} \frac{1}{2r} \int_{-r}^r \| \Psi(t) \|_\alpha
\ dt &\leq \lim_{r \to \infty} \frac{M(\alpha)}{2r}\int_{-r}^r
\int_{0}^{+\infty}s^{-\alpha} e^{-\gamma s} \ \| \zeta(t-s) \| ds
\ dt\\
&\leq \lim_{r \to \infty} M(\alpha)\int_{0}^{+\infty} s^{-\alpha}
e^{-\gamma s} \frac{1}{2r}\int_{-r}^r \| \zeta (t -s) \| dt \
ds=0,
\end{align*}
by using Lebesgue dominated Convergence theorem, and the fact that
$AP_0(\X)$ is invariant under translations. Thus  $\Psi$ belongs
to $AP_0(\X_\alpha)$.

The proof for $\Gamma_4u(\cdot)$ is similar to that of $\Gamma_3
u(\cdot)$. However one makes use of \eqref{hyQ*} rather than
\eqref{hyP*}.
\endproof

\begin{lemma}\label{ll3} Under assumptions {\rm ({\bf H1})-({\bf H2})-({\bf H3})-({\bf
H4})}, the integral operators $\Gamma_1$ and $\Gamma_2$ defined
above map $PAP(\X_\alpha)$ into itself.
\end{lemma}

\proof Let $u \in PAP(\X_\alpha)$. Since $B \in B(\X_\alpha,\X)$
it follows that the function $t \mapsto B u(t)$ belongs to
$PAP(\X)$. Again, using the composition theorem of pseudo almost
periodic functions \cite[Theorem 5]{am} it follows that $\psi
(\cdot) = f (\cdot , Bu(\cdot))$ is in $PAP(\X_\beta)$ whenever $u
\in PAP(\X_\alpha)$. In particular,
$\displaystyle{\|\psi\|_{\infty,\beta} = \sup_{t \in \R} \| f( t ,
Bu(t))\|_{\beta} < \infty}.$

 Now write $\psi = w + z,$ where
$w \in AP(\X_\beta)$ and $z \in AP_0(\X_\beta)$, that is,
$\Gamma_1 \phi = \Xi (w) + \Xi (z)$ where
$$
\Xi w(t) := \int_{-\infty}^{t} A T(t -s) Pw(s) ds, \ \ \mbox{and}
\ \ \Xi z(t) := \int_{-\infty}^{t} A T(t -s) Pz(s) ds.
$$

Clearly, $\Xi(w) \in AP(\X_\alpha)$. Indeed, since $w \in
AP(\X_\beta)$, for every $\varepsilon > 0$ there exists
$l(\varepsilon)
> 0$ such that for all $\xi$ there is $\tau \in [\xi , \xi
+ l(\varepsilon)]$ with the property: $$\|w(t+\tau) - w(t)
\|_{\beta} < \nu\varepsilon \ \ \mbox{for each} \ t \in \R,$$
where $\nu=\dfrac1{ c\gamma^{\beta-\alpha}\Gamma(\beta-\alpha)}$
and $c$ being the constant appearing in Lemma \ref{beta}.

Now, the estimate \eqref{beta2} yields
\begin{align*}
\left\| \Xi(w)(t + \tau) - \Xi(w)(t)\right\|_\alpha& =\left\|\int_0^{+\infty} A T(s) P\left(w(t-s + \tau) - w(t-s)\right)\,ds\right\|_\alpha\\
&\leq\int_0^{+\infty} s^{\beta-\alpha-1}e^{-\gamma s} \|w(t-s + \tau) - w(t-s)\|_\beta ds\\
&\leq \varepsilon
\end{align*}
for each $t \in \R$, and hence $\Xi(w) \in AP(\X_\alpha)$.

Now, let $r > 0$. Then, by \eqref{beta2}, we have
\begin{eqnarray*}
\frac{1}{2r} \int_{-r}^r \|(\Xi z)(t)\|_{\X_\alpha} \,dt &\leq&
\frac{1}{2r}\int_{-r}^{r}\int_{0}^{+\infty} \|AT(s) Pz(t-s)
\|_{\alpha} \,ds \,dt \\
\\&\leq& \frac{1}{2r}\int_{-r}^{r}\int_0^{+\infty} s^{\beta-\alpha-1}e^{-\gamma s}
\|z(t-s)\|_{\beta} \, ds\,dt \\&\leq&
\int_0^{+\infty}s^{\beta-\alpha-1}e^{-\gamma
s}\frac{1}{2r}\int_{-r}^{r} \|z(t-s)\|_{\beta} \,dt \, ds.
\end{eqnarray*}

Obviously, $\displaystyle \lim_{r \to \infty} \frac{1}{2r}
\int_{-r}^r \|(\Xi z)(t)\|_{\alpha} dt = 0,$ since $t \mapsto
z(t-s) \in AP_0(\X_\beta)$ for every $s\in \R$. Thus $\Xi z \in
AP_0(\X_\alpha)$.

The proof for $\Gamma_2u(\cdot)$ is similar to that of $\Gamma_1
u(\cdot)$. However, one uses \eqref{beta1} instead of
\eqref{beta2}.
\endproof

Throughout the rest of the paper, the constant $k(\alpha)$ denotes
the bound of the embedding $\X_\beta \hookrightarrow \X_\alpha$,
that is,
$$\|u\|_\alpha \leq k(\alpha) \|u\|_\beta \ \ \mbox{for each} \ u \in
\X_\beta.$$

\begin{theorem}\label{theo}Under the  assumptions {\bf (H1)-(H2)-(H3)-(H4)}, the evolution equation \eqref{0} has
a unique pseudo almost periodic mild solution whenever $\Theta<1$,
where
\[
\Theta =  K \varpi \left[ k(\alpha) + \frac{c}{\delta} +
c\frac{\Gamma(\beta-\alpha)}{\gamma^{\beta-\alpha}}+
\frac{M(\alpha) \Gamma(1-\alpha)} { \gamma^{1-\alpha}} +
\frac{c(\alpha) } { \delta} \right],
\]
and $\varpi = \max(\|B\|_{B(\X_\alpha,\X)}, \|C\|_{B(\X_\alpha,
\X)})$.
\end{theorem}

\proof Consider the nonlinear operator ${\mathbb M}$ on
$PAP(\X_\alpha)$ given by

\begin{eqnarray}
{\mathbb M} u(t)&=&-f(t,Bu(t)) - \int_{-\infty}^{t}AT(t-s)P f(s,
Bu(s))ds
\nonumber\\
&+& \int_{t}^{\infty}AT(t-s)Qf(s, Bu(s))ds +
\int_{-\infty}^{t}T(t-s)Pg(s, Cu(s))ds \nonumber\\
&-& \int_{t}^{\infty}T(t-s)Qg(s, Cu(s))ds\nonumber
\end{eqnarray}
for each $t \in \R$.

As we have previously seen, for every $u\in PAP(\X_\alpha)$,
$f(\cdot,Bu(\cdot))\in PAP(\X_\beta)\subset PAP(\X_\alpha)$. In
view of Lemma \ref{ll2} and Lemma \ref{ll3}, it follows that
${\mathbb M}$ maps $PAP(\X_\alpha)$ into itself. To complete the
proof one has to show that ${\mathbb M}$ has a unique fixed-point.

Let $v,w\in PAP(\X_\alpha)$
\begin{eqnarray*}
\left\| \Gamma_1 (v)(t) - \Gamma_1(w)(t)\right\|_\alpha &\leq&
\int_{-\infty}^t
\left\|AT(t-s)P\left[f(s,Bv(s))-f(s,Bw(s))\right]\right\|_\alpha
ds
\nonumber \\
&\leq& cK \|B\|_{B(\X_\alpha,\X)} \|v - w\|_{\infty, \alpha}
\int_{-\infty}^t (t-s)^{\beta-\alpha-1}e^{-\gamma(t-s)} ds \\
&=&  c\frac{\Gamma(\beta-\alpha)}{\gamma^{\beta-\alpha}}K
\|B\|_{B(\X_\alpha,\X)}\|v - w\|_{\infty, \alpha}.
\end{eqnarray*}

Similarly,

\begin{eqnarray*}
\| \Gamma_2 (v)(t) - \Gamma_2(w)(t)\|_\alpha &\leq&
\int_{t}^\infty
\|AT(t-s)Q\left[f(s,Bv(s))-f(s,Bw(s))\right]\|_\alpha ds
\nonumber \\
&\leq& cK \|B\|_{B(\X_\alpha,\X)} \|v - w\|_{\infty, \alpha}
\int_{t}^{+\infty} e^{\delta(t-s)} ds\nonumber \\
&=& \frac{cK \|B\|_{B(\X_\alpha,\X)}}{\delta} \|v - w\|_{\infty,
\alpha}.
\end{eqnarray*}

Now for $\Gamma_3$ and $\Gamma_4$, we have the following
approximations
\begin{eqnarray}
\| \Gamma_3 (v)(t) - \Gamma_3(w)(t)\|_\alpha &\leq&
\int_{-\infty}^t
\|T(t-s)P\left[g(s,Cv(s))-g(s,Cw(s))\right]\|_\alpha ds
\nonumber \\
&\leq& \frac{K \|C\|_{B(\X_\alpha, \X)} M(\alpha)
\Gamma(1-\alpha)} { \gamma^{1-\alpha}} \|v - w\|_{\infty, \alpha},
\nonumber
\end{eqnarray}
and

\begin{eqnarray}
\| \Gamma_4 (v)(t) - \Gamma_4(w)(t)\|_\alpha &\leq&
\int_{t}^{+\infty}
\|T(t-s)Q\left[g(s,Cv(s))-g(s,Cw(s))\right]\|_\alpha ds
\nonumber \\
&\leq& K c(\alpha) \|C\|_{B(\X_\alpha, \X)} \|v - w\|_{\infty,
\alpha} \int_{t}^{+\infty} e^{\delta(t-s)} ds\nonumber \\
&=& \frac{K  \|C\|_{B(\X_\alpha, \X)} c(\alpha) } { \delta}\|v -
w\|_{\infty, \alpha}. \nonumber
\end{eqnarray}

Consequently,
$$
    \|\mathbb{M}v-\mathbb{M}w\|_{\infty, \alpha}\leq \Theta \,.\,\|v -
    w\|_{\infty, \alpha}.$$
Clearly, if $\Theta < 1,$ then (\ref{0}) has a unique fixed-point
by the Banach fixed point theorem, which obviously is the only
pseudo almost periodic solution to (\ref{0}).
\endproof

\begin{example}
For $\sigma \in \R$, consider the (semilinear) transport equation
with delay endowed with Dirichlet conditions:

\begin{eqnarray}\label{toka1}
\hspace{1cm} \displaystyle \frac{\partial}{\partial
t}[\varphi(t,x) + f(t, \varphi(t-p,x))]&=&
\frac{\partial}{\partial x}\varphi(t,x)+ \sigma \varphi(t,x) +
g(t,\varphi(t - p,x))\\
\displaystyle \varphi(t,0)&=& \varphi(t, 1) = 0, \ \ t \in
\R\label{toka2}
\end{eqnarray}
where $p > 0$, and $f, g: \R \times C[0 , 1] \mapsto C[0 , 1]$ are
some jointly continuous functions.

Set $\X:=C([0,1])$ equipped with the sup norm. Define the operator
$A$ by
$$A(\varphi):= \varphi' + \sigma \varphi, \ \ \forall \varphi \in D(A),$$ where $ D(A):=\{ \varphi\in C^1([0,1]),
\varphi(0)= \varphi(1)= 0\} \subset C[0 , 1]$.

 Clearly $A$ is sectorial, and hence is the generator of an analytic
semigroup.  In addition to the above, the resolvent and spectrum
of $A$ are respectively given by
$$\rho(A) = \C - \{ 2n \pi i + \sigma: \ n \in \Z\} \ \ \mbox{and}
\ \ \sigma(A)= \{ 2 n \pi i + \sigma: \ n \in \Z\}$$ so that
$\sigma(A) \cap i \R = \{\emptyset\}$ whenever $\sigma \not= 0$.

\begin{theorem}Under assumptions {\bf (H3)-(H4)}, if $\sigma \not= 0$, then the
transport equation with delay \eqref{toka1}-\eqref{toka2} has a
unique $\X_{\alpha}$-valued pseudo almost periodic mild solution
whenever $K$ is small enough.

\end{theorem}

\end{example}

\begin{example}
For $\sigma \in \R$, consider the (semilinear) heat equation with
delay endowed with Dirichlet conditions:

\begin{eqnarray}\label{toka4}
\hspace{1cm} \displaystyle \frac{\partial}{\partial
t}[\varphi(t,x) + f(t, \varphi(t-p, x))]&=&
\frac{\partial^2}{\partial x^2}\varphi(t,x)+ \sigma \varphi(t,x) +
g(t,\varphi(t-p,x))\\
\displaystyle \varphi(t,0)&=& \varphi(t, 1) = 0\label{toka5}
\end{eqnarray}
for $t \in \R$ and $x \in [0, 1]$, where $p > 0$, and $f,g : \R
\times C[0 , 1]\mapsto C[0 , 1]$ are some jointly continuous
functions.

Take
 $\X:=C[0 , 1]$, equipped with the sup norm. Define the operator
$A$ by
$$A(\varphi):= \varphi'' + \sigma\varphi, \ \ \forall \varphi \in D(A),$$ where $ D(A):=\{ \varphi\in C^2[0,1],
\varphi(0)= \varphi(1)= 0\} \subset C[0 , 1]$.

 Clearly $A$ is sectorial, and hence is the generator of an analytic
semigroup.  In addition to the above, the resolvent and spectrum
of $A$ are respectively given by
$$\rho(A) = \C - \{ - n^2 \pi^2 + \sigma: n \in \N \} \ \ \mbox{and}
\ \ \sigma(A)= \{ - n^2 \pi^2 + \sigma: n \in \N \}$$ so that
$\sigma(A) \cap i \R = \{\emptyset\}$ whenever $\sigma \not= n^2
\pi^2$.

\begin{theorem}Under assumptions {\bf (H3)-(H4)}, if $\sigma \not= n^2 \pi^2$ for each $n \in \N$, then the
heat equation with delay \eqref{toka4}-\eqref{toka5} has a unique
$\X_{\alpha}$-valued pseudo almost periodic mild solution whenever
$K$ is small enough.

\end{theorem}

\end{example}

\bigskip

{\bf Acknowledgement.} The author wants to express many thanks to
Professor Maniar for useful comments and suggestions on the
manuscript.

\bigskip

\bibliographystyle{amsplain}

\end{document}